\documentclass[a4paper,11pt]{amsart}
\usepackage[utf8]{inputenc}
\usepackage[english]{babel}
\usepackage{amsthm}
\usepackage{amsfonts}
\usepackage{amsmath}
\usepackage{amssymb}
\calclayout

\usepackage{graphicx}

\newtheorem{theorem}{Theorem}

\newcommand{\Q}{\mathbb{Q}}

\newcommand{\A}{\mathcal{A}}

\begin{document}

\title{Rational Right Triangles of a Given Area}
\markright{Rational Right Triangles of a Given Area}
\author{Stephanie Chan}
\address{Department of Mathematics, University College London, Gower Street, London, WC1E~6BT, United Kingdom}
\email{stephanie.chan.16@ucl.ac.uk}

\begin{abstract}
Starting from \textsl{any} given rational-sided, right triangle, for example the $(3,4,5)$-triangle with area $6$, we use Euclidean geometry to show that there are infinitely many other rational-sided, right triangles of the same area. We show further that the set of all such triangles of a given area is finitely generated under our geometric construction. Such areas are known as ``congruent numbers'' and have a rich history in which all the results in this article have been proved and far more. Yet, as far as we can tell, this seems to be the first exploration using this kind of geometric technique.
\end{abstract}
\maketitle

\section{Introduction}
\begin{figure}[ht]
\centering
\includegraphics{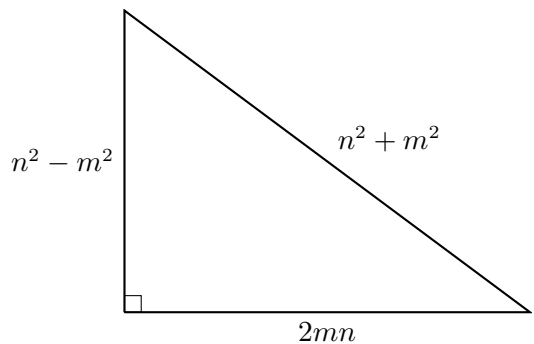}
\caption{Euclid's parameterisation of integer-sided right triangles.}
\label{fig:euclid}
\end{figure}

Book 10, Proposition 29 of Euclid's \textsl{Elements} \cite{Euclid} gives the familiar formula for finding right triangles with integer sides:\ For any integers $n>m\geq 1$ we have the right triangle shown in Figure~\ref{fig:euclid}.

\noindent Every integer-sided right triangle can be found from this formula and subsequently multiplying
through by an integer scalar (for example, we obtain the $(9,12,15)$-triangle by taking $m=1,n=2$ and tripling each side). Moreover \emph{all rational-sided} right triangles can be obtained by scaling these triangles by rational multiples.

Five hundred years after Euclid, Diophantus noted that $X,Y,Z$ is an integer-sided right triangle of (integer) area $\A$ if and only if $Z^2-4\A$ and $Z^2+4\A$ are both squares \cite[Chapter XVI]{Dickson}. Viewed like this, it is of interest to determine whether there are any solutions for a given integer $\A$, and if so, to find all the triples 
$X,Y,Z$ with a given area $\A$. This was one of the earliest questions in mathematical research, worked on not only by the ancient Greeks, but also by Arab mathematicians of the tenth century (including al-Karaji), by Leonardo Fibonacci in thirteenth century Pisa, and by many others ever since.

The $(3,4,5)$-triangle is the only integer-sided right triangle of area $6$. Multiplying through by $70$ we obtain the $(210,280,350)$-triangle which has the same area as the $(49,1200,1201)$-triangle. This begs the question as to whether there are areas $\A$ for which there are arbitrarily many integer-sided right triangles of area $\A$. To attack this question it is easiest to work with rational-sided right triangles, so that the 
$(\frac{49}{70},\frac{1200}{70},\frac{1201}{70})$-triangle also has area $6$. We say that $\A$ is a \emph{congruent number} if it is the area of a rational-sided right triangle. Fibonacci's first work on this question was to answer John of Palermo's 1220 challenge, namely to find a rational square $w^2$ such that $w^2-5$ and $w^2+5$ are both squares; he found $w=\frac {41}{12}$ which corresponds to the triangle
$(\frac 32,\frac{20}3,\frac{41}6)$ of area $5$, and this scales up to the familiar right triangle $(9,40,41)$ of area $5\cdot 6^2$.

The first few Pythagorean triples give rise to the congruent numbers
\[
 5, 6,7, 13, 14, 15, 20, 21, 22, 24, 28, 30, 34, 39, 41, 46, \ldots.
\]
It is known that $23, 29, 31, 37, 38, 45, 47,\ldots$ are also congruent, and each number involves triangles with surprisingly large numerators and denominators in their side lengths. For example, the smallest numbers in a rational right triangle of area $23$ are 
 $\left(\frac{80155}{20748}, \frac{41496}{3485}, \frac{905141617}{72306780}\right)$.

Given a congruent number, like $6$, we might ask to find all integer right triangles that have area of the form
$6r^2$. We already saw $(3,4,5)$ and $(49,1200,1201)$, and the next example is $(2896804,7216803,7776485)$, followed by an even larger triangle $(339252715200, 2066690884801, 2094350404801)$. Maybe we can keep finding more and more such triples, but are there infinitely many of them? The integers involved seem to get larger rapidly, so that a brute force search is likely to become very inefficient. Is there a better way of generating such triples?

There is a simple connection between congruent numbers and the theory of elliptic curves:\ Suppose that the $(2mn,n^2-m^2,n^2+m^2)$-triangle, scaled by a factor $r$, has area $\A$. Then 
\[
y^2 = x(x^2-\A^2), 
\]
where $x=\A n/m$ and $y=\A^2/rm^2$ are rationals. That is, we have obtained a rational point on this \emph{elliptic curve}, the set of points $(x,y)$ on $y^2 = x(x^2-\A^2)$, which we denote by $E_\A$. Modern research into congruent numbers typically studies this \emph{congruent number curve} with many of the tools of arithmetic geometry
(see, e.g., \cite{Koblitz}).
This connection with elliptic curves gives answers to our questions, beautifully but indirectly. 
We will obtain some of the results with Greek geometry, in particular given a rational-sided triangle of area $\A$, we will construct a different rational-sided triangle of the same area, with bigger numerator and denominator. This process can be iterated to give an infinite sequence of such triangles, and we can show that the triangles are all nonsimilar.
We will follow up later by relating this process to the theory of elliptic curves.

Elliptic curve theory implies that right-angled triangles must exist with many of the properties that we determine in this article, and even gives their side lengths. The theory does not indicate how to construct these new triangles from the old, however. Our more direct approach constructs these triangles explicitly. This allows us to visualise the proofs through these triangles, giving an alternative view, which we hope has the potential to make the relevant concepts more tractable.

\subsection*{Acknowledgement}
The author would like to thank her Ph.D. advisor, Andrew Granville, for his guidance and inspiring discussions in the development of this paper. She is also grateful to Jennifer Balakrishnan, John Coates and Henri Darmon for their helpful comments and suggestions.

\section{Geometric development of the parameterisation of a right triangle}\label{section:par}
\begin{figure}[ht]
\centering
\includegraphics[scale=0.95]{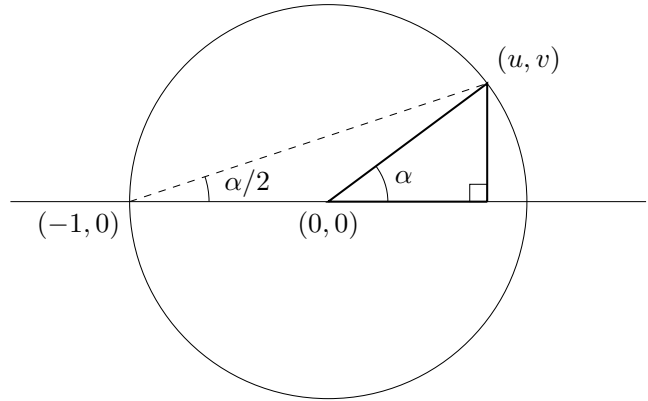}
\caption{Parameterisation of a right triangle.}
\label{fig:par}
\end{figure}

Given a right triangle with rational sides $X,Y,Z$ and angle $\alpha$ opposite the side of length $Y$, we rescale the triangle by setting $u=X/Z$ and $v=Y/Z$, so that the hypotenuse has length $1$, and 
 $u=\cos\alpha$ and $v=\sin\alpha$ are rational numbers. We embed the triangle inside a circle of radius $1$, so that the vertex with angle $\alpha$ lies at the origin, and the hypotenuse is a radius of the circle, as in Figure~\ref{fig:par}.
 
\noindent We extend the line on the $x$-axis to the end of the circle and then by classic Greek geometry this vertex subtends an angle of $\alpha/2$. The number $t=\tan(\alpha/2)$ is the slope of the line from $(-1,0)$ to $(u,v)$, and so 
\[
t=\frac v{u+1} \in (0,1),
\]
and is therefore also a rational number.
The double angle formulae express trigonometric functions at $2x$ as rational functions of $\tan x$:
\[\sin 2x=\frac{2\tan x}{1+\tan^2 x},\quad \cos 2x=\frac{1-\tan^2 x}{1+\tan^2 x},\quad \tan 2x =\frac{2\tan x}{1-\tan^2 x}.\]
Therefore we get the natural parameterisation,
\[u=\frac{1-t^2}{1+t^2},\quad v=\frac{2t}{1+t^2}.\]
Note that we can obtain the same parameterisation entirely algebraically, without referencing the trigonometric functions.
Writing the rational number $t$ as $m/n$ where $m$ and $n$ are coprime positive integers, and scaling up the triangle by a factor $m^2+n^2$ yields Euclid's parametrization given in Figure~\ref{fig:euclid}.

If we switch the roles of $X$ and $Y$, then the triangle would be parametrized with the value of $t$ replaced by $(1-t)/(1+t)$. Since the map $t\mapsto (1-t)/(1+t)$ is its own inverse, there are exactly two parameters in the range $(0,1)$ for each triangle. 
However, there is a primitive choice for each triangle, which is the unique parameter with numerator and denominator of different parities. Indeed, given a triangle, write the two parameters as $t=m/n$ and $(1-t)/(1+t)=M/N=(n-m)/(n+m)$ with $(m,n)=(M,N)=1$. The only possible common factor of $n-m$ and $n+m$ is $2$.
If $m$ and $n$ are of different parities, then $M+N=(n-m)+(n+m)=2n$ is even.
If $m$ and $n$ are both odd, then $M+N=(n-m)/2+(n+m)/2=n$ is odd.
Therefore, in exactly one of $t$ and $(1-t)/(1+t)$, the numerator and denominator have different parities. Such parameters correspond to primitive triangles in Figure~\ref{fig:euclid} (otherwise all three sides are divisible by $2$).

\section{Obtaining a new triangle geometrically}\label{section:double}
Given a rational right triangle $ABC$ of given area $\A$, our aim is to geometrically construct a genuinely different rational right triangle $A'B'C'$ with area $r^2\A$ for some rational $r$.

\begin{figure}[ht]
\centering
\includegraphics{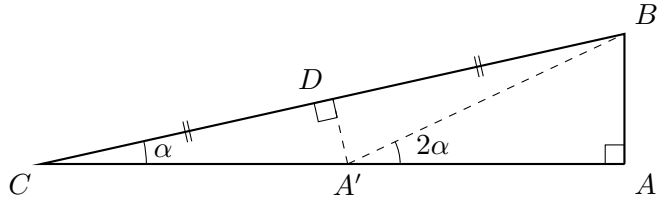}
\caption{Doubling the angle $\alpha$ of the initial triangle $ABC$.}
\label{fig:double1}
\end{figure}

As above, let $\alpha$ be the acute angle of $ABC$, so that $\sin \alpha, \cos \alpha$ and $\tan \alpha$ are all rational numbers. The perpendicular bisector of the hypotenuse $BC$ must intersect the longer leg of the triangle at some point $A'$ as shown in Figure~\ref{fig:double1}.

\noindent The geometry implies that $A'C=A'B$, and that the angle $BA'A$ is double the angle $BCA$, and so equals $2\alpha$. Moreover $\sin 2\alpha$, $\cos 2\alpha$ and $\tan 2\alpha$ are all rational expressions of $\tan\alpha$ and so are themselves rational.
Since $AB$ is rational we deduce that the triangle $A'AB$ has rational sides, and so $A'C=AC-A'A$ is rational. 

\begin{figure}[ht]
\centering
\includegraphics{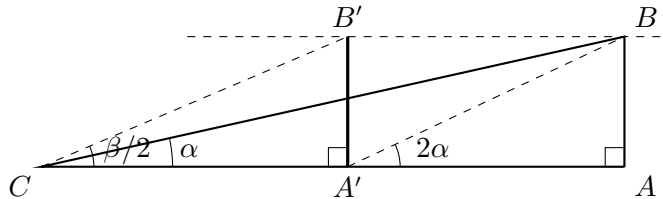}
\caption{Obtaining the parameter $\tan(\beta/2)=\sin 2\alpha$ of the new triangle.}
\label{fig:double2}
\end{figure}

We next translate the side $AB$ to $A'B'$, keeping it perpendicular to $AC$, and let $\beta/2$ denote the angle $B'CA'$ (See Figure~\ref{fig:double2}). 

\noindent Therefore,
\[
\tan{\left(\frac{\beta}{2}\right)} = \frac{A'B'}{A'C}=\frac{AB}{A'B}=\sin 2\alpha,
\]
which we have already noted is rational.

\begin{figure}[ht]
\centering
\includegraphics{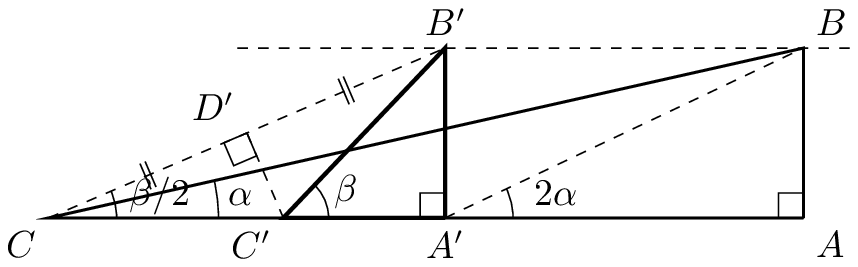}
\caption{Completing the new triangle $A'B'C'$ with angle $\beta$.}
\label{fig:double3}
\end{figure}

Since $A'C=A'B>AB=A'B'$, we see that $A'C$ is the longer leg of the right triangle $A'B'C$.
Therefore, the perpendicular bisector of the hypotenuse $B'C$ must intersect the longer leg of the triangle, $A'C$, at some point $C'$ as illustrated in Figure~\ref{fig:double3}.

\noindent Now $B'C'A'$ has angle $\beta$, by Euclidean geometry. 
As $\tan(\beta/2)$ is rational, we deduce that $\sin \beta, \cos \beta$ and $\tan \beta$ are all rational
by the double angle formulae. Moreover, since $A'B'=AB$ is rational, we deduce that all the sides of the triangle $A'B'C'$ are rational. 

The legs $A'B'$ of the right triangle $A'B'C'$, and $AB$ of $ABC$, have equal length, and so
the ratio of the areas of the triangles $A'B'C'$ and $ABC$ is
\[\begin{split}
\frac{ A'C'}{ AC}
&=\frac{AB/AC}{A'B'/A'C'}
=\frac{\tan\alpha}{\tan\beta}
=\frac{1-\tan^2(\beta/2)}{2\tan(\beta/2)}\cdot\frac{\sin\alpha\cos\alpha}{\cos^2\alpha}\\
&=\frac{1-\sin^2 2\alpha}{2\sin 2\alpha}\cdot\frac{\sin 2\alpha}{2\cos^2\alpha}
=\frac{\cos^2 2\alpha}{4\cos^2\alpha}
={\left(\frac{\cos 2\alpha}{2\cos \alpha}\right)}^2={\left(\frac{A'A\cdot BC}{2\cdot A'B \cdot AC}\right)}^2,
\end{split}\]
a rational square, as claimed.

In Section~\ref{section:double} we saw that the parameter for the triangle $ABC$ is $t=\tan(\alpha/2)$; and therefore the 
parameter of the new triangle $A'B'C'$, is 
\[
T=\tan{\left(\frac{\beta}{2}\right)}= \sin 2\alpha=2\sin\alpha\cos\alpha=2uv=\frac{4t(1-t^2)}{(1+t^2)^2} \in (0,1) 
\]
by the double angle formulae. If we replace $t$ by its associated parameter $\frac{1-t}{1+t}$ in this formula, we obtain the same value of $T$. 

\section{A little algebra} 
\label{section:algebra}
\subsection{Iterating the construction} 

For a given rational right triangle of area $\A$, we determine the associated parameter $t$, as in Section~\ref{section:par}. Write $t_1=t$, and then construct a new rational right triangle of area $\A r_2^2$ for some rational $r_2$ with parameter $t_2=T$ by the method of Section \ref{section:double}. Then repeat this construction and create an infinite sequence of rational numbers $t_1,t_2,\ldots,$ where $t_k$ is the 
 parameter for a 
rational right triangle of area $\A r_k^2$ for some rational $r_k$, and
$t_{k+1}=4t_k(1-t_k^2)/(1+t_k^2)^2$.

If we start with the triple $(3, 4, 5)$, then the next triple is $(49, 1200, 1201)$, followed by $(339252715200, 2066690884801, 2094350404801)$. The numbers involved grow quickly. Motivated by this observation, we will show that the parameters of the triangles have strictly increasing denominators, and therefore iterating the construction cannot produce any similar triangles.

We write the rational number $t_k=m_k/n_k$, where $m_k$ and $n_k$ are coprime positive integers. To establish that the $t_j$ are all distinct, and therefore give rise to different triangles, we will prove that $n_1<n_2<n_3<\cdots$.

Now given $t=m/n$ with $(m,n)=1$ where $m$ and $n$ have different parities, then $T=M/N$ where $M=4mn(n^2-m^2)$ and $N=(m^2+n^2)^2$ are coprime, and $M$ is even and $N$ is odd. This implies that, for all $k\geq 1$, $n_{k+1}=(m_k^2+n_k^2)^2$ and so
$n_k^4< n_{k+1}<4n_k^4$ (as $m_k<n_k$ since $t_k\in (0,1)$). We have therefore proved the result we were aiming for:

\begin{theorem} \label{infinite}
If there is one rational-sided right triangle of area $\A$ then there are infinitely many.
\end{theorem}

\subsection{Descent} The integer-sided triangles constructed from the $(3,4,5)$-triangle, by the method given in Section~\ref{section:double}, have legs $7^2$ and $3\cdot 20^2$, and then $3\cdot 336280^2$ and $1437599^2$, respectively. In each example one leg is a square, the other is $2\A$ times a square. This always happens in this construction, for if we began with a right triangle with coprime integer sides $X,Y,Z$, then the new triangle with coprime integer sides constructed in Section~\ref{section:double}, has legs of length 
\[
 (X^2-Y^2)^2 \ \ \text{and} \ \ 4XYZ^2=2\A \cdot (2Z)^2.
\]
These expressions equal $(1-T^2)\cdot Z^4$ and $2T\cdot Z^4$, respectively, so that $1-T^2$ is a square, and even more is true:
\[
1-T= 1-2uv=(u-v)^2\ \text{ and } \ 1+T= 1+2uv=(u+v)^2 
\]
are both squares of rational numbers. So the new triangle, constructed in Section~\ref{section:double}, has a very particular arithmetic form, which suggests that perhaps this is the only way such a triangle can arise.

\begin{theorem} [Descent] \label{doubleweil}
Given a right triangle with sides $(1-T^2,2T,1+T^2)$ for which $1-T$ and $1+T$ are both rational squares, there exists a right triangle with sides $(u,v,1)$ such that $T=2uv$, where $u$ and $v$ are both rational. The ratio of their areas is the square of a rational number.
\end{theorem}

We call this a \emph{descent} since the numerator and denominators in the $(u,v,1)$-triangle are smaller than those in the first triangle, as we proved in the previous subsection.

\begin{proof} 
Write $1-T=r^2$ and $1+T=s^2$, where $r$ and $s$ are rational.
Then let $u=(r+s)/2$ and $v=(s-r)/2$, so that $u^2+v^2=\frac 12(r^2+s^2)=1$ and 
$2uv = {\frac 12 (s^2-r^2)}=T$, as desired. 

Given the $(u,v,1)$-triangle we can reconstruct the original triangle by the method of Section~\ref{section:double}, and so the ratio of their areas is a square.
\end{proof}

In every newly constructed triangle of this sort, we see the corresponding integer-side right triangle has a square leg: Write $T=M/N$, and then one leg has length $N^2-M^2=N^2(1-T)(1+T)$ which is a square by
Theorem \ref{doubleweil}. However, there are integer-side right triangles that do not arise in this way, but which have a square leg, for example $(9,40,41)$: Here $t=4/5$ so that although $1-t^2=(3/5)^2$, we have
$1-t=1/5$ and $1+t=9/5$, neither of which is a square, so Theorem \ref{doubleweil} does not apply.

\section{The group law on an elliptic curve}
\begin{figure}[ht]
\centering
\includegraphics[scale=0.8]{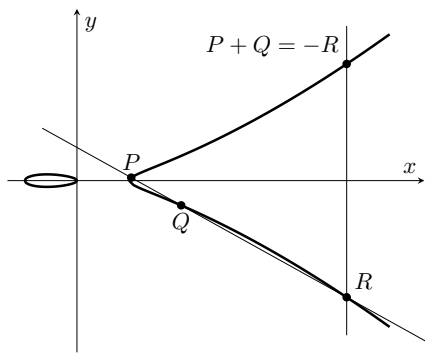}
\caption{Group law on the elliptic curve $y^2=x^3-6^2x$.}
\label{fig:group}
\end{figure}

Any equation of the form $y^2=x^3+ax+b$, for which $4a^3+27b^2\ne 0$ defines an \emph{elliptic curve}, which we denote by $E$. The rational points $(x,y)$ together with $\mathcal{O}$, the point at infinity, are denoted by $E(\mathbb Q)$. Poincar\'e showed that the complex points on $E$ form a group, and that $E(\mathbb Q)$ is a subgroup. The group law, a construction that was discovered much earlier, is quite extraordinary \cite[Chapter 1]{SilvermanTate}:
Given $P,Q\in E(\mathbb Q)$ we draw the line between them, whose equation must have rational coefficients, as the points $P$ and $Q$ are rational. Any line intersects a cubic curve in three points, so suppose the third point is $R$ and the line has equation $y=mx+c$. Then the $x$-coordinates of $P,Q$ and $R$ must all satisfy the equation
\[
x^3-(mx+c)^2+ax+b = x^3+ax+b -y^2 = 0.
\]
There are three roots to this equation, $x(P), x(Q)$ and $x(R)$, and their sum is $m^2$, and so $x(R)$ is rational. Finally as $R$ lies on the line $y=mx+c$, we deduce that $R \in E(\mathbb Q)$. The point $P+Q$ is then given by reflecting $R$ in the $x$-axis, as in Figure~\ref{fig:group}.

\noindent There is one special case to consider more carefully: If $P=Q$, then we let our line be the tangent at $P$ (which is what would happen if we brought $Q$ close in towards $P$). It is still true that the line intersects the curve at three points, but now we must make sure to count those points \emph{with multiplicity}, so we say that the line intersects the curve twice at $P$. We construct $2P=P+Q$ by again taking the third point of intersection $R$, and reflecting it in the $x$-axis.

For example, the tangent line at $P=(12,36)\in E_6(\mathbb Q)$ to the curve $E_6$ is $y=\frac {11}2x-30$ so that
$x(R)+2x(P)=\frac {121}4$ and therefore $x(R)=\frac {121}4-24=\frac {25}4$. Thus $y(R)=\frac {11}2\cdot \frac {25}4-30=\frac{35}8$ and therefore $2P=(\frac {25}4,-\frac{35}8)$.

\subsection{Using the congruent number curve to find our new triangle}\label{section:map} The set of rational right triangles with area $\A$ can be represented, as we saw in Section~\ref{section:par}, by the set 
\[
T_\A:=\{ t\in \mathbb Q:\ t(1-t^2) = \A r^2 \ \text{for some} \ r \in \mathbb Q\}.
\]
We define a map $\Phi: T_\A \to E_\A(\mathbb Q)$ by $\Phi(t)=(\A/t,r\A^2/t^2)$. This map is easily inverted since $\Phi^{-1}(x,y)=\A/x$. We can make sense of any negative elements $-1/t\in(-\infty,-1)$ in this set by identifying each with $t\in(0,1)$. In this way the negative part $T_\A\cap(-\infty,-1)$ can be viewed as merely a copy of $T_\A\cap(0,1)$. 

\begin{figure}[ht]
\centering
\includegraphics{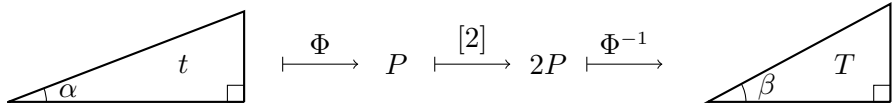}
\caption{new congruent triangles via the multiplication-by-2 map.}
\label{fig:doublemap}
\end{figure}

One way to find a new triangle in $T_\A$ from a given $t\in T_\A$ is to first determine $P=\Phi(t)$, then compute $2P$, and finally let $T=\Phi^{-1}(2P)\in T_\A$; see Figure~\ref{fig:doublemap}.
 
For example, starting with the $(3,4,5)$ triangle, we have $t=\frac 12\in T_6$ and then $P=\Phi(t)=(12,36)\in E_6(\mathbb Q)$. At the end of the last subsection we found that $2P=(\frac {25}4,-\frac{35}8)$, and so $T=\Phi^{-1}(2P)=\frac {24}{25} \in T_6$, the same $T$-value we obtained from our geometric construction!

In general we have $x(P)=\A/t$ and then $x(2P)=(\frac{1+t^2}{2r})^2$, so that $T=\frac{4\A r^2}{(1+t^2)^2}=
\frac{4t(1-t^2)}{(1+t^2)^2}$. This explains in terms of arithmetic geometry the results we obtained above.

\section{Everything in twos}
\subsection{The two theorems in the context of elliptic curves}

For any point $P\in E_\A (\mathbb Q)$ with $x(P)\ne 0,\A ,-\A $ one can obtain $t=\Phi^{-1}(P)\in T_\A $ and Theorem \ref{infinite} then implies that the infinite sequence of points $t_1=t, t_2, \ldots \in T_\A $ are all distinct. Defining $P_n=\Phi(t_n)$, we obtain an infinite sequence of distinct points $P_1=P,P_2, \ldots \in E_\A (\mathbb Q)$ where $P_n=2^nP$, so $P$ has infinite order inside the group of rational points $E_\A (\mathbb Q)$. This implies that the only points of finite order on $E_\A (\mathbb Q)$ are the point $\mathcal O$ and the three points $(0,0), (\A ,0), (-\A ,0)$ of order two.

Theorem \ref{doubleweil} can be re-interpreted as stating that $P\in 2E_\A (\mathbb Q)$ if and only if $1-T$ and $1+T$ are both squares, where $T=\Phi^{-1}(P)$, and shows us how to determine a point $R=\Phi(\frac v{u+1})\in E_\A (\mathbb Q)$ for which $P=2R$. Now any solution to $2Q=P$ must be of the form $Q=R+S$ where $2S=\mathcal O$, and so $Q=R, R+(0,0), R+(\A ,0)$, or $R+(-\A ,0)$. These points correspond to taking the four solutions $(\pm u,\pm v,1)$ in Theorem \ref{doubleweil}, but to correspond to an actual triangle we restrict to the one case where $u$ and $v$ are both positive.

Fermat proved that $1$ is not a congruent number using a method called infinite descent, which we can think of as applying Thereom \ref{doubleweil} repeatedly to get smaller and smaller triangles, eventually reaching a contradiction. For an exposition on Fermat's method, see \cite[Chapter II, \S X]{Weil}.
Inspired by Fermat's method, Mordell proved that for any elliptic curve $E$, the group $E(\mathbb Q)$ is finitely generated (see \cite[Chapter 16]{Mordell}). His proof was technically difficult, but this situation was improved a few years later by Weil, who introduced several important ideas to better understand and simplify Mordell's proof. More about that later.

\subsection{Two different maps}
\begin{figure}[ht]
\centering
\includegraphics{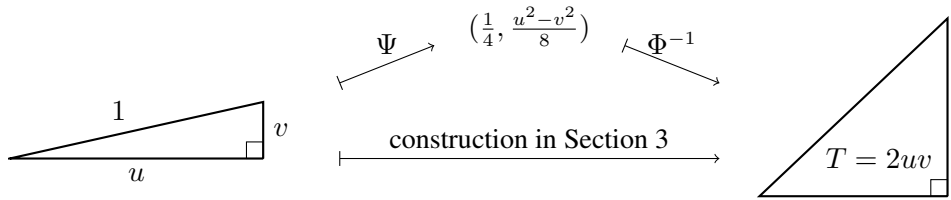}
\caption{Composing two different maps.}
\label{fig:compose}
\end{figure}

The original approach of Diophantus shows that $(X,Y,Z)$ is a rational right triangle of area $\A $ if and only if
\[
Z^2+4\A =(X+Y)^2 \ \text{ and } \ Z^2-4\A =(X-Y)^2.
\]
Multiplying these together and dividing by $16$ gives 
\[ r^4-\A ^2=s^2 \ \text{ where } \ r=Z/2 \ \text{ and } \ s=(X^2-Y^2)/4.\]
This then yields the rational point $(r^2,rs)$ on the congruent number curve 
$E_\A :y^2=x^3-\A ^2x$. 
Call this Diophantus-inspired map that takes $(X, Y, Z)$ to $(r^2, rs)$ $\Psi$.
In Section~\ref{section:map}, we presented a different map $\Phi$ using the parameter $t$, yielding the rational point $(\A /t,r\A ^2/t^2)$ on $E_\A $. This map $\Phi$ is easily inverted giving $t=\A /x$.
So let's combine these two maps: We start with a triangle with parameter $t$ and construct $(u,v,1)$ of area $\A $. The map $\Psi$ yields the rational point $(1/4,(u^2-v^2)/8)$ on $E_\A $, and inverting the second map yields $T=4\A =2uv$ as before; see Figure~\ref{fig:compose}.

\subsection{Two different $t$-values for the same triangle} 
We saw that every right triangle is parameterised by two possible values in $(0,1)$, namely $t$ and $\psi(t):= \frac{1-t}{1+t}$. 
Now $\psi(t)(1-\psi(t)^2)=\A(\frac{2r}{(1+t)^2})^2$ and so 
$\Phi(\psi(t))=(\frac{\A (1+t)}{1-t}, 2r (\frac{\A }{1-t})^2)$. 
The map  $\Phi(t)\to \Phi(\psi(t))$ is a special case of the map $(x,y)\to(\frac{\A (x+\A )}{x-\A }, 2y( \frac{ \A }{x-\A })^2)$, which is an \emph{isogeny} of order two, and commutes with the multiplication-by-2 map.

\section{Adding two different points}
\begin{figure}[ht]
\centering
\includegraphics{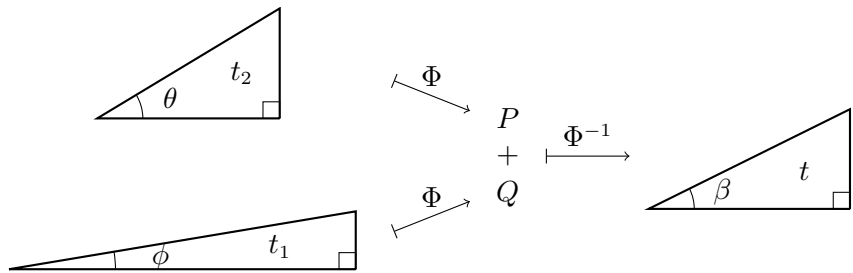}
\caption{Adding points corresponding to nonsimilar triangles on the congruent number curve.}
\label{fig:add}
\end{figure}

Adding two different points $P$ and $Q$ on $E_\A (\mathbb Q)$ should correspond to somehow combining two rational right-sided triangles of area $\A $, to create another rational right-sided triangle of area $\A $. We now present a geometric construction to do this directly, as suggested by Figure~\ref{fig:add}.

\begin{figure}[ht]
\centering
\includegraphics{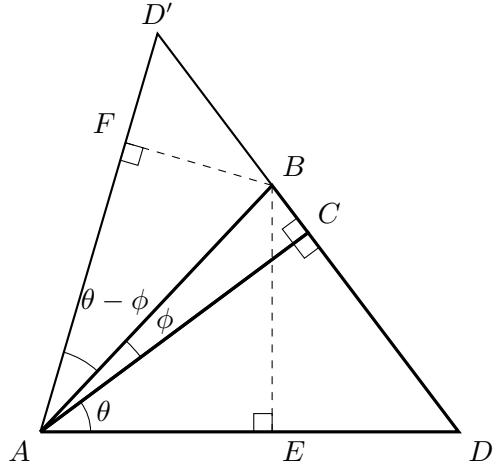}
\caption{Obtaining the angles $\theta\pm\phi$.}
\label{fig:add1}
\end{figure}

We begin with two nonsimilar rational right triangles of the same area. We rescale the initial triangles so that their longer legs have the same length, and then align them as in Figure~\ref{fig:add1}: The two triangles are $ABC$ with angle $\phi$, and $ACD$ with angle $\theta$, where $AB=1$. Next we reflect the point $D$ through the line $AC$ to obtain the point $D'$, and the angles $\theta-\phi$ and $\theta+\phi$. See Figure~\ref{fig:add1}.

\noindent As $AB=1$, we deduce that $\sin \phi =BC,\ \cos \phi=AC$ and so 
$CD = \tan \theta \cos \phi$, $AD=\frac{\cos \phi}{\cos \theta}$, which are all rational numbers.
By the formulae for adding and subtracting angles, we also deduce that the trigonometric functions $\sin, \cos$ 
and $\tan$, evaluated at $\theta-\phi$ and $\theta+\phi$, are rational numbers.

The triangle $ABC$ has area $\frac 12\cdot AC \cdot BC = \frac 12 \sin \phi \cos \phi = \frac 14 \sin 2\phi$, and the triangle $ACD$ has area $\frac 12\cdot AC\cdot CD=\frac 12 \cos^2 \phi \tan \theta = \frac 14 (\frac{\cos \phi}{\cos \theta})^2 \sin 2\theta$. Our hypothesis states that the ratio of these areas is the square of a rational number, and so their product must be a square, which implies that $\sin 2\phi\sin 2\theta$ is the square of a rational, say $R^2$.

\begin{figure}[ht]
\centering
\includegraphics{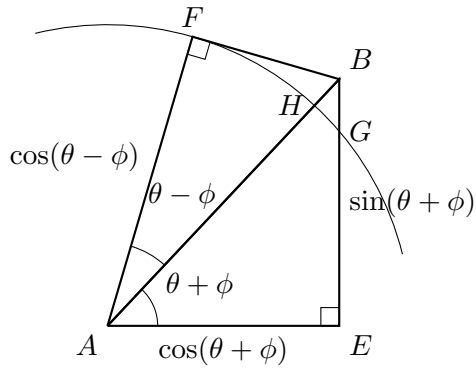}
\caption{Constructing the lengths $AH=\cos(\theta-\phi)$ and $EG=R$.}
\label{fig:add2}
\end{figure}

We drop a perpendicular from $B$ to the lines $AD$ and $AD'$,
and then focus on the triangles $ABE$ and $ABF$. The circle centred at $A$ with radius $AF$ intersects the line $AB$ at $H$, and $BE$ at $G$; see Figure~\ref{fig:add2}.

\noindent Now $AH=AG=AF=\cos(\theta-\phi)$ and $BH=AB-AH=1-\cos(\theta-\phi)$.
Moreover, \[EG=\sqrt{AG^2-AE^2}=\sqrt{\cos^2(\theta-\phi)-\cos^2(\theta+\phi)}=\sqrt{\sin 2\theta \sin 2\phi}=R,\] so that $BG=BE-GE=\sin(\theta+\phi)-R$.

\begin{figure}[ht]
\centering
\includegraphics{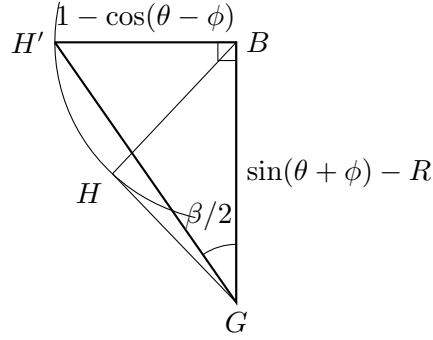}
\caption{Rotating the side $BH$ to get a right triangle.}
\label{fig:add3}
\end{figure}

The triangle $BGH$ contains the information we need, specifically, the lengths $BH$ and $BG$. Rotate the side $BH$ to $BH'$, so that $BH'$ is perpendicular to $BG$; see Figure~\ref{fig:add3}.

\noindent
Let $\beta/2$ be the angle at the corner $BGH'$, so that 
\[
\tan{\left(\frac{\beta}{2}\right)}=\frac{1-\cos(\theta-\phi)}{\sin(\theta+\phi)-R} ,
\] 
and is therefore a rational number.

\begin{figure}[ht]
\centering
\includegraphics{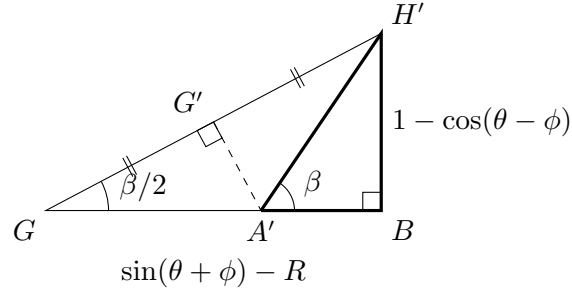}
\caption{Constructing the triangle from parameter $\tan(\beta/2)$.}
\label{fig:add4}
\end{figure}

The perpendicular bisector of $GH'$ intersects $BG$ at a point $A'$, yielding a new triangle $A'BH'$ with angle $\beta$ as shown in Figure~\ref{fig:add4}.

\noindent
This new triangle $A'BH'$ is rational since it has a rational side ${BH'}{=}{1-\cos(\theta-\phi)}$, and
$A'B=BH'/\tan \beta$ and $A'H'=BH'/\sin \beta$, where $\sin \beta$ and $\tan \beta$ are rational functions of 
$\tan(\beta/2)$ and so are themselves rational. The parameter for the triangle $A'BH'$ is
\begin{equation}\label{eq:addition}
\begin{split}
t&=\tan{\left(\frac{\beta}{2}\right)}
=\frac{ (t_1-t_2)^2}{t_2(1-t_1^2)+t_1(1-t_2^2)-2\sqrt{ t_1 t_2 (1-t_1^2) (1-t_2^2) }} \\
& =\left( \frac{ t_1-t_2}{ \sqrt{t_2(1-t_1^2)} -\sqrt{t_1(1-t_2^2) }} \right)^2
=\left( \frac{ \sqrt{t_2(1-t_1^2)} +\sqrt{t_1(1-t_2^2) }}{ 1+ t_1t_2} \right)^2,
\end{split}
\end{equation} 
where $t_1$ is the parameter for $ABC$ and $t_2$ is the parameter for $ADC$ 
so, up to a square, they have areas $T_1=\frac{4t_1(1-t_1^2)}{(1+t_1^2)^2}=\sin 2\phi$ and
$T_2=\frac{4t_2(1-t_2^2)}{(1+t_2^2)^2}=\sin 2\theta$ so that $R^2=T_1T_2$.
In both $t_1$ and $t_2$ the numerator and denominator have different parities. It is evident, by writing out the two expressions for $t$, in terms of the numerators and denominators of $t_1$ and $t_2$, that the numerator and denominator of $t$ also have different parities.

The area of the new triangle is
\[\begin{split}
\frac{4\cdot\text{Area}(A'BH')}{|A'H'|^2} &=T=\frac{4t(1-t^2)}{(1+t^2)^2} \\
&={\left(\frac{{\left(4t_1t_2-(1-t_1^2)(1-t_2^2)\right)}{\left(\sqrt{T_1} +\sqrt{T_2}\right)}}{4t_1t_2+(1-t_1^2)(1-t_2^2)+2(1-t_1t_2)(t_1+t_2)\sqrt{T_1T_2}}\right)}^2.
\end{split}\]
Therefore, if $T_1=\A r_1^2$ and $T_2=\A r_2^2$, then $T=\A r^2$ where
\[
 r= \frac{{\left(4t_1t_2-(1-t_1^2)(1-t_2^2)\right)}(r_1 +r_2)}{4t_1t_2+(1-t_1^2)(1-t_2^2)+2(1-t_1t_2)(t_1+t_2)\A r_1r_2}.
\]

\subsection{Subtraction} 
We obtain a different triangle if we subtract a point $P$ from a different point $Q$ on $E_\A(\Q)$. Two different triangles with parameters $t_1$ and $t_2$ will result in a new triangle with parameter 
\begin{equation}\label{eq:subtraction}
t=\left( \frac{ \sqrt{t_2(1-t_1^2)} -\sqrt{t_1(1-t_2^2) }}{ 1+ t_1t_2} \right)^2,
\end{equation}
of area $T=\frac{4t(1-t^2)}{(1+t^2)^2}=\A r^2$ where
\[
 r= \frac{{\left(4t_1t_2-(1-t_1^2)(1-t_2^2)\right)}(r_1 -r_2)}{4t_1t_2+(1-t_1^2)(1-t_2^2)-2(1-t_1t_2)(t_1+t_2)\A r_1r_2}.
\]
The same parameter is obtained by first adding the triangles with respect to the parameters $\psi(t_1)=\frac{1-t_1}{1+t_1}$ and $t_2$, that is, taking $\theta$ to be the other angle of the triangle $ACD$ as in the geometric construction outlined before, then applying $\psi$ to the resulting parameter.
This process is also a geometric construction.
We call this subtraction because if we are given parameters $t_1$ and $t_2$ that yield $t_3$ under \eqref{eq:addition}, we can recover $t_1$ by applying \eqref{eq:subtraction} to $t_3$ and $t_2$.

\section{The set of rational right triangles with area $\A $}
Let $\mathcal{T}_\A $ denote the set of all rational right triangles of area $\A $.
This set is in 1-to-1 correspondence with the set of parameters $t=m/n\in T_\A$, where $m$ and $n$ are coprime positive integers of different parities. We define addition and subtraction on $\mathcal{T}_\A $ as described geometrically in the previous section so that the new triangles have parameters given by \eqref{eq:addition} and \eqref{eq:subtraction} respectively.
One can show that if we begin with parameters with numerator and denominators of different parities, the new parameter obtained from \eqref{eq:addition} or \eqref{eq:subtraction} has the same property. 
We will henceforth restrict our attention to such parameters in $T_\A $.
One has to be a little careful in that when carrying out subtractions on triangles, $\triangle_1-\triangle_2$ and $\triangle_2-\triangle_1$ are geometrically the same triangle.

\subsection{A finite set of generators} 


We define a map $W:T_\A \to\Q^*/(\Q^*)^2\times\Q^*/(\Q^*)^2$ by
$W(t)=(1-t,1+t)$. Two rational numbers are equal in $\Q^*/(\Q^*)^2$ if their ratio is a rational square.
This map is a homomorphism, for if $\triangle_1+\triangle_2=\triangle_3$ and correspond to parameters $t_1,t_2$ and $t_3$ respectively, we obtain from the expression for $t_3$ in \eqref{eq:addition},
\[\begin{split}
(1+t_1)(1+t_2)(1+t_3)
&=\left(\frac{(1+t_1)(1+t_2)+\sqrt{t_1t_2(1-t_1^2)(1-t_2^2)}}{1+ t_1t_2}\right)^2\\
&=\left((1+t_1)(1+t_2)+\frac{1}{4}\A r_1r_2(1+t_1^2)(1+t_2^2)\right)^2
=1
\end{split}\]
in $\Q^*/(\Q^*)^2$. Similarly, we have $(1-t_1)(1-t_2)(1-t_3)=1$ in $\Q^*/(\Q^*)^2$. 

Theorem \ref{doubleweil} implies that $2T_\A $ is the kernel of $W$; that is, $t\in 2T_\A $ if and only if $W(t)=(1,1)$. Moreover the equivalence classes $T_\A /2T_\A $ are given by $\{ t\in T_\A :\ W(t)=w\}$ for each $w\in W(T_\A )$, the image of $T_\A $ under the map $W$.

Now if $t=m/n$ with $(m,n)=1$, then $mn(m+n)(n-m)=n^4t(1-t^2)=\A (rn^2)^2$. The four factors
$m,n,m+n,n-m$ are pairwise coprime. Therefore $m,n,m+n,n-m=a,b,c,d$ in $\Q^*/(\Q^*)^2$ where $abcd=\A $ in $\Q^*/(\Q^*)^2$, and so
$W(t)=(bd,bc)$ for some coprime, squarefree integers $b,c,d$ where $bcd$ is a divisor of $\A $.
This restricts $W(T_\A )$ to a finite set of possibilities and so implies that $T_\A /2T_\A $ is finite.

Pick a complete set $\mathcal{R}$ of representatives of equivalence classes in $T_\A /2T_\A $.
Fix any $t$ in $T_\A $; then $t$ must be in the same equivalence class with some $T$ in $\mathcal{R}$. Denote the corresponding triangles of $t$ and $T$ by $\triangle$ and $\triangle'$, respectively. Now $\triangle+\triangle'$ is in the kernel of $W$, so there exists some $\triangle_1$ with parameter $t_1$ in $T_\A $ such that
\[\triangle+\triangle'=2\triangle_1.\]
Given $t=m/n$, $t_1=m_1/n_1$ and $T=M/N$ with $(m,n)=(m_1,n_1)=(M,N)=1$, the triangle $\triangle+\triangle'$ has parameter \[\left(\frac{\sqrt{NM(n^2-m^2)}+\sqrt{nm(N^2-M^2)}}{nN+mM}\right)^2.\]
The denominator divides and hence is smaller than $(nN+mM)^2<4n^2N^2$.
On the other hand, we also know that the parameter of $2\triangle_1$ has denominator greater than $n_1^4$ from Section~\ref{section:algebra}. Combining the two inequalities gives us $n_1^4<4n^2N^2$. 

Take $C$ to be the maximum of the denominators of all elements in $\mathcal{R}$. Then if $n\geq 2C$, we have the inequality $n_1<\sqrt{2nN}\leq\sqrt{2Cn}\leq n$.
We repeat this process on $t_1$ and inductively obtain a sequence of parameters $t_2,t_3,\dots$ defining triangles $\triangle_2,\triangle_3,\dots$, respectively, so that for each $k\geq 1$,
\[\triangle_{k}+\triangle'_{k}=2\triangle_{k+1},\]
where $\triangle'_{k}$ has parameter in $\mathcal{R}$ and $t_k$ has denominator $n_k$.
Provided $n_{k}\geq 2C$, the denominators form a strictly decreasing sequence of integers $n_{k+1}<n_{k}<\dots<n_1<n$. This must yield a $t_K$ with denominator $n_K<2C$ after a finite number of steps and we obtain a linear relation
\[\triangle+\triangle'+2\triangle'_{1}+\dots+2^{K-1}\triangle'_{K-1}=2^K\triangle_{K}.\]
Parameters in $T_\A $ with denominators less than $2C$ form a finite set containing $\mathcal{R}$ and generate all the triangles in $\mathcal{T}_\A $. 

We have therefore proved the following result.

\begin{theorem}
There exist triangles $\triangle_1,\dots,\triangle_r$ in $\mathcal{T}_\A $ such that for any triangle ${\triangle\in\mathcal{T}_\A}$, there exists integers $a_1,\dots,a_r$ for which
\[
\triangle = a_1\triangle_1+\cdots +a_r\triangle_r .
\]
\end{theorem}

What we have shown is essentially Mordell's theorem, but in the setting of our geometrically constructed set of triangles. The original version of Mordell's theorem states that $E_\A (\Q)$ is a finitely generated abelian group of rank $r$, say. That is, there is a basis $P_1,\dots,P_r$ for $E_\A (\Q)$, and this must correspond to some basis of triangles $\triangle_1,\dots,\triangle_r$ for $\mathcal{T}_\A $.
Our map $W$ is a special case of the Weil map on $E_\A :y^2=x^3-\A ^2x$, which is the homomorphism 
\[\begin{split}
&E_\A (\Q)\to\Q^*/(\Q^*)^2\times\Q^*/(\Q^*)^2\times\Q^*/(\Q^*)^2\\
&(x,y)\mapsto(x,x-\A ,x+\A ),\end{split}\]
and is well known to have kernel $2E_\A (\Q)$ (see \cite[Chapter X.1]{Silverman}). 
The Weil map is inspired by Fermat's method of descent.
It plays a role analogous to our map $W$ in performing descent in the proof of Mordell's theorem.

It is an open question as to how we might find the set $\mathcal{R}$. Although there is a finite explicit set of possibilities for $W(T_\A )$ and thus for $T_\A /2T_\A $, one does not know, after a failed finite search for an example in a given coset, whether one has failed because there is no example, or because one has not searched far enough. This is an example of one of the outstanding problems in computational arithmetic geometry, determining the rank of the group $E(\Q)$ of a given elliptic curve $E$. In practice, it is usually possible to compute the rank of a given explicit example, but no known algorithm is guaranteed to work in general.

\subsection{Partitioning in terms of parameter $t$} 


We note that 
since $\triangle_1,\dots,\triangle_r$ are independent, this implies that 
\[
 a_1\triangle_1+\cdots +a_r\triangle_r:\ 0\leq a_1,\ldots,a_r\leq 1
 \]
 gives a complete set of representatives of $\mathcal{T}_\A /2\mathcal{T}_\A $ and $|\mathcal{T}_\A /2\mathcal{T}_\A |=2^r$.

A proportion of $1/2^r$ of the triangles have that both $1-t$ and $1+t$ are squares, and so
at least $1/2^r$ of the triangles have $1-t^2$ equal to a square. This corresponds, in Figure~\ref{fig:euclid}, to
$m^2-n^2=m^2(1-t^2)$ being a square. That is, at least $1/2^r$ of the primitive triangles have one leg having square length.
 
\subsection{First example}
We first look at the case when $\A =6$.
Now $E_6(\Q)$ is generated by $P=(12,36)$, which corresponds to the triangle $\triangle_P$ with sides $(3,4,5)$ and parameter $\frac 12$.
The even multiples of $\triangle_P$ have $(1-t,1+t)=(1,1)\mod(\Q^*)^2$ and the odd multiples of $\triangle_P$ have $(1-t,1+t)=(2,6)\mod(\Q^*)^2$. 
The set $\mathcal{R}$ is $\{0,\frac 12\}$, so $C=2$. There are only $3$ possible fractions with denominator less than $2C=4$: $\frac 12,\frac 13,\frac 23$ and $\triangle_P$ is the only triangle from these parameters that is in $\mathcal{T}_6$. 

\subsection{Second example}
The area $\A =34$ is the smallest integer such that $r>1$. The group $E_{34}(\Q)$ has $r=2$, generated by the points $P=(\frac{289}4,\frac{4335}8)$ and $Q=
(578,13872)$, which corresponds to the triangles $\triangle_P=(225, 272, 353)$ and $\triangle_Q=(17, 144, 145)$, respectively.
The triangles in $\mathcal{T}_{34}$ can be partitioned into four sets with their representatives tabulated in the following table. 
\begin{center}
 \begin{tabular}{ | c | c | c | c | c | c | }
 \hline
 Representatives & $t$ & {$W(t)$} \\
 \hline
 $0$ & $0$ & $(1, 1)$ \\ \hline
 $\triangle_P=(225, 272, 353)$ & $\frac{8}{17}$ & $(17, 17)$ \\ \hline
 $\triangle_Q=(17, 144, 145)$ & $\frac{8}{9}$ & $(1,17)$ \\ \hline
 $\triangle_P+\triangle_Q=(1377,3136, 3425)$ & $\frac{32}{49}$ & $(17, 1)$ \\ \hline
 \end{tabular}
\end{center}
The set $\mathcal{R}$ is $\{0,\frac{8}{17},\frac{8}{9},\frac{32}{49}\}$ and $C=49$.
By checking all possible fractions, we find that $\triangle_P$, $\triangle_Q$, $\triangle_P+\triangle_Q$ are the only triangles in $\mathcal{T}_{34}$ which have parameters with denominators less than $2C=98$.

\section{Futher thoughts}
There are other scenarios in which our geometric approach can be applied.

Adding a number in the quadratic field $\Q(\sqrt{d})$ to its conjugate gives a rational number. 
Similarly, adding a point on $E_\A (\Q(\sqrt{d}))$ to its conjugate gives a point on $E_\A (\Q)$. This gives rise to a rational right triangle as long as the new point is nontrivial. Our geometric construction allows us to construct this new triangle directly. In this way we can establish the following.

\begin{theorem} \label{addconj}
If there exists a right triangle with side lengths in $\Q(\sqrt{d})$ of rational area $\A $ such that the rational part of its hypotenuse is nonzero, then there is a rational-sided right triangle of area $\A $.
\end{theorem}

Subtracting a point on $E_\A (\Q(\sqrt{d}))$ from its conjugate gives a new point $P$ on $E_{\A} (\Q(\sqrt{d}))$ with conjugate $-P$. This can be reinterpreted as a point on $E_{\A d} (\Q)$ giving rise to a rational-sided triangle of area $\A d$. Again, this new triangle can be constructed using our geometric construction, giving the following result.

\begin{theorem} \label{subconj}
If there exists a right triangle with side lengths in $\Q(\sqrt{d})$ and irrational hypotenuse of rational area $\A $, then there is a rational-sided right triangle of area $\A d$.
\end{theorem}

Are there other applications of our geometric method in understanding ideas related to the congruent number problem?
Can it give an alternative insight into some of the other intriguing theorems and conjectures discussed in \cite{Koblitz}? We hope the interested reader will join us on this voyage of discovery.


\begin{thebibliography}{1}

\bibitem{Dickson}
L.~E.~Dickson,
\textit{History of the Theory of Numbers}. Vol.~{II}, Diophantine Analysis,
Chelsea Publishing Co., New York, 1966.

\bibitem{Euclid}
Euclid,
\textit{Euclid's Elements: All Thirteen Books Complete in One Volume}.
Trans. by T.~L.~Heath.
Ed. by D.~Densmore.
Green Lion Press, Santa Fe, NM, 2002.

\bibitem{Koblitz}
N.~Koblitz,
\textit{Introduction to Elliptic Curves and Modular Forms}.
Graduate Texts in Mathematics, Vol.~97, Springer-Verlag, New York, 1984.

\bibitem{Mordell}
L.~J.~Mordell,
\textit{Diophantine Equations}.
Pure and Applied Mathematics, Vol.~30, Academic Press, London-New York, 1969.

\bibitem{Silverman}
J.~H.~Silverman,
\textit{The Arithmetic of Elliptic Curves}.
Second edition.
Graduate Texts in Mathematics, Vol.~106, Springer, Dordrecht, 2009.

\bibitem{SilvermanTate}
J.~H.~Silverman, J.~T.~Tate,
\textit{Rational Points on Elliptic Curves}.
Second edition.
Undergraduate Texts in Mathematics, Springer, Cham, 2015.

\bibitem{Weil}
A.~Weil,
\textit{Number Theory: An Approach through History from Hammurapi to Legendre}.
Reprint of the 1984 edition. 
Modern Birkh\"auser Classics, Birkh\"auser Boston, Inc., Boston, MA, 2007. 

\end{thebibliography}
\end{document}